\documentclass[12pt]{article}
\usepackage{amsmath,amssymb}
\newtheorem{tm}{Theorem}[section]
\newtheorem{lm}[tm]{Lemma}

\newtheorem{re}[tm]{Remark}

\newtheorem{pr}[tm]{Proposition}
 
 \newenvironment{demo}[1]{\par\smallskip\par\begin{trivlist}
\item[]{\bf #1}\ }{\end{trivlist}\par\smallskip\par}
\newcommand{\Proof}{\begin{demo}{{\it Proof.\ }}}
\newcommand{\qed}{\end{demo}}
\newcommand{\toy}{\ \rule[0em]{0.5ex}{1.8ex}}
\newcommand{\QED}{\toy\end{demo}}

\newcommand{\la}{\langle}
\newcommand{\ra}{\rangle}
\newcommand{\nn}{\nonumber}
\newcommand{\III}{{\vert \kern-.10em \vert \kern-.10em \vert}}
\newcommand{\ve}{\varepsilon}

\makeatletter
 
 \@addtoreset{equation}{section}
\makeatother
 
\begin{document}
\setlength{\baselineskip}{15pt} 
%
\bibliographystyle{plain}
\title{
A moment estimate of the derivative process in rough path theory
}
\author{{\Large Yuzuru INAHAMA
}\\
Graduate School of Mathematics,   Nagoya University\\
Furocho, Chikusa-ku, Nagoya 464-8602, JAPAN 
\\
E-mail:~\tt{inahama@math.nagoya-u.ac.jp}
}
\date{   }
%
%
\maketitle

\begin{center}
{\bf Abstract}
\end{center}
In this paper we prove the derivative process of a  rough differential equation
driven by Brownian rough path
has finite $L^r$-moment for any $r \ge 1$.
Thanks to Burkholder-Davis-Gundy's inequality,  this kind of problem is easy in the usual SDE theory. 
In the context of rough path theory, however, it does not seem so obvious.

\vspace{5mm}
%
%
\section{Introduction }
In stochastic analysis,    the derivative process of a given stochastic differential equation 
(or equivalently, equation of the stochastic flow)
has been studied extensively, 
because  it plays a  very important role  in various situations.
On the other hand, in rough path theory, the derivative process was not studied very much.
One reason could be that it has unbounded coefficients, for  existence of solution for 
a rough differential equation (RDE) is in general difficult.
The aim of this paper is to prove $L^r$-integrability for the first level path of the derivative process for any $r \ge 1$.

Now we will give a more detailed explanation.
We consider the following  RDE in a Banach setting. Here, $W$ is  Brownian rough path and $a$ is an initial value.
\[
dY_t  =  \sigma (Y_t)  dW_t,  \qquad\qquad  Y_0=a. 
\]
Its derivative equation is given by
\[
dJ_t  =  \nabla \sigma (Y_t)  \la J_t  \bullet , dW_t \ra,  \qquad\qquad  J_0={\rm Id}
\]
Roughly speaking, $J_t$ is the derivative of a map $a \mapsto Y_t =Y(a)_t$.
It is known among experts of this research field that a unique solution $(Y,J)$ in the rough sense exists,
although  there seems to be no published paper which proves it.

\vspace{5mm}

Our main result is as follows (see Theorem \ref{tm.main} below for details);
Let $2<p<3$ and let $W$ be  Banach space-valued Brownian rough path. 
Then, under a suitable condition on the coefficient $\sigma$,
$1/p$-H\"older norm of the first level path $J^1$ is $L^r$-integrable for any $r \ge 1$.

\vspace{5mm}

This kind moment estimate appears in various occasions.
In the author's case, when he and H. Kawabi 
try to prove a stationary phase for solutions RDEs in a forthcoming paper \cite{ik} 
(which can be regarded as  a rough path version of \cite{ba}), 
this type of integrability of the derivative process $J$ is needed.

To the author's knowledge, the only exposition which explicitly  discusses the derivative equation is Aida's 
unpublished manuscript \cite{aida}, in which he established T. Lyons' continuity theorem for $(Y,J)$
and proved the following estimate; If $W$ is controlled by a control function $\omega$, i.e., $|W^i_{s,t}| \le \omega(s,t)^{i/p}$
for all $0 \le s \le t \le 1$ and $i =1,2$,
then $\sup_t |J^1_{0,t}| \le  C \exp (C \omega (0,1)^{\gamma})$ for some constants $C>0$ and $\gamma \ge 1$.
This kind of deterministic argument is of great importance.
However, it is  not suitable for our purpose, becsuse even if we take $\omega(s,t)= \|  W^1 \|^{p}_{p-var; [s,t]} + \|  W^2 \|^{p/2}_{p/2-var; [s,t]} $,
the right hand side is not integrable. 
Here, $\|  \,\cdot\,\|_{p-var; [s,t]} $ is $p$-variation norm on the subinterval $[s,t]$.
Therefore, we need to take a closer look at the argument in \cite{aida}.

\begin{re}
(i)~
the same results holds for an RDE with a drift term. 
In such a case,  we only need to consider the lift of a "space-time" process  $t \mapsto (w_t, t)$  instead of Brownian rough path.
Here, $w$ is Brownian motion (i.e., the first level path of $W$).
\\
\noindent
(ii)~
The author does not know whether the main result  is true or not when the driving rough path is a lift of fractional Brownian motion.
See Lemma \ref{lm.mom.Mall} and remark \ref{re.fracok} below for details.
\end{re}

%
%
\section{Setting}
Let  $( {\cal V}, {\cal H}, \mu)$ be an abstract Wiener space and 
let $(w_t)_{0 \le t \le 1}$ be Brownian motion on ${\cal V}$ associated to $\mu$, which starts at $0$.
Let $2 < p <3$ and let $G\Omega_p ({\cal V})$ be the geometric rough path space over ${\cal V}$ with $p$-variation norm. 
(When given the $1/p$-H\"older norm, the geometric rough path space is denoted by $G\Omega_{1/p - Hld} ({\cal V})$.) 
In this article,  the time interval is always $[0,1]$ and tensor spaces of Banach spaces are equipped with the projective
tensor norm, unless stated otherwise.
We basically use the original formulation of rough path theory as in Lyons and Qian \cite{lq},
although there are a few variants of the theory now.

We denote by $w(m)$ be the $m$th dyadic approximation of $w$, i.e., 
the piecewise linear approximation associated to the partition $\{ 0 <1/2^m < 2/2^m <\cdots <  (2^m-1)/2^m <1  \}$.
Its lift, i.e., the smooth rough path above $w(m)$, is denoted by $W(m)$ as usual.
Unlike the finite dimensional case, existence of Brownian rough path is not known.
So we set the following assumption:
\\
\\
\noindent
{\bf (A1)} \qquad
 $\lim_{n,m \to \infty} {\mathbb E}[ \| W(m)^i  - W(n)^i  \|_{i/p-Hld}   ] $ 
for  $i=1,2$.  
Here, the norm is  H\"older norm of  index $i/p$.
\\
\\
The limit is denoted by $W$ and is called Brownian rough path.
(This formulation is used in Deriech \cite{der}.  The well-known sufficient condition "Exactness " as in Definition 4.6.1, \cite{lq} 
implies {\bf (A1)}.
So, we will work under this assumption.)
Interestingly, {\bf (A1)} implies almost sure convergence of  $W(m)$, too (see \cite{der}).

%
%
%
Let ${\cal W}$ be another real Banach space and $\sigma: {\cal W} \to L({\cal V} , {\cal W} )$
be $C^4_b$ in the Fr\'echet sense (i.e., $\nabla^j \sigma$ is bounded for $j=0,1, \ldots ,4$).
Here, $L({\cal V} , {\cal W} )$ denotes the set of bounded linear maps from ${\cal V}$ to ${\cal W}$, 
which is equipped with the operator norm.
Consider the following RDE;  for $X \in G\Omega_p ({\cal V})$,
\begin{equation}\label{eq.rde.ori}
dY_t =  \sigma(Y_t) dX_t,  \qquad \mbox{ with $Y_0 = a \in {\cal W}$}
\end{equation}
A solution in the rough path sense is $Z =(X,Y) \in G\Omega_p ({\cal V} \oplus {\cal W})$.
Its second component  $Y \in G\Omega_p ({\cal W} )$ is also called a solution.
Under this regularity condition for $\sigma$, this RDE has a unique solution.
So, $X \mapsto Y (= :\Phi(X))$ defines a map, which is called It\^o map. 
By T. Lyons' continuity theorem, $\Phi : G\Omega_p ({\cal V} ) \to G\Omega_p ({\cal W} )$ 
is (locally Lipschitz) continuous.

Aida \cite{aida} gave a rahter quantative estimate for the growth of
 the solution and the local  Lipschitz constant for two solutions.
If $X$ is controlled by a control function $\omega$  
(that is,  $|X^i_{s,t} | \le \omega(s,t)^{i/p}$ for all $s \le t$ and $i=1,2$),  
then  $Z=(X,Y)$ is controlled by a control function $\hat{\omega}$ of the form
$\hat{\omega}(s,t) =C(1+\omega(0,1)^{\gamma}) \omega(s,t)$ with certain positive constants $C, \gamma$, 
which is independent of the initial value $a$.

Adding to RDE (\ref{eq.rde}),
we also consider the following ``derivative equation.''
\begin{equation}\label{eq.rde.der}
dJ_t =  \nabla \sigma(Y_t) \la J_t \bullet , dX_t \ra
 \qquad \mbox{with $J_0 ={\rm Id}_{{\cal W}},$}
\end{equation} 
Notice that  a solution $J$ takes its values in $L({\cal W}):=L({\cal W}, {\cal W})$
and that
$\nabla \sigma(Y_t) $ is a bounded bilinear map from ${\cal W} \times {\cal V}$ to ${\cal W}$.
Formally,  $J_t$ is the derivative of a map $a \in {\cal W} \mapsto a + Y^1_{0,t}  \in {\cal W}$.

%
%
RDEs (\ref{eq.rde.ori}) and  (\ref{eq.rde.der}) are  obviously equivalent to consider the following RDE;  
\begin{eqnarray}
dY_t &=&  \sigma(a + Y_t) dX_t,  \qquad \mbox{ with $Y_0 =0 \in {\cal W}$}
\label{eq.rde}\\
dJ_t &=&  \nabla \sigma(a+Y_t) \la ( {\rm Id}_{{\cal W}} +  J_t )\bullet , dX_t \ra
 \qquad \mbox{with $J_0 ={\rm Id}_{{\cal W}},$}
 \label{eq.rde.der.sh}
  \end{eqnarray}
We will basically study RDEs  (\ref{eq.rde}) and  (\ref{eq.rde.der.sh}).
(More precisely, the shitfed equations as above are the definition of RDE with a non-zero initial condition.)

RDE (\ref{eq.rde.der.sh}) can be written in a simpler way as follows;
\begin{equation}\label{eq.rde.mat.sh}
dJ_t =   d   (  {\rm Id}_{{\cal W}} + J_t )  = (dM_t) \cdot (  {\rm Id}_{{\cal W}} + J_t )  \qquad \mbox{with $J_0 =0,$}
\end{equation}%
where $M$ an $L({\cal W})$-valued path, which is given  by 
\begin{equation}\label{eq.rde.M}
M_t := \int_0^t \nabla \sigma(a +Y_u) \la \,\cdot\,, dX_u \ra.
\end{equation}
Note that the map  $\Lambda : {\cal V} \oplus {\cal W}  \to L ({\cal V} \oplus {\cal W}  , L({\cal W})  )$
defined by
\[
\Lambda (x,y)  \la x' , y' \ra  =
\nabla \sigma(y ) \la \,\cdot\,, x' \ra    \quad  \in L({\cal W})
 \]
is $C^3_b$ and the right hand side of (\ref{eq.rde.mat.sh}) is well-defined as  rough path integral.
So,  if $X$ is controlled by $\omega$, then 
the right hand side of (\ref{eq.rde.M}) is a rough integral 
and $M$ is controlled by $\omega '(s,t) =C'(1+\omega(0,1)^{\gamma '}) \omega(s,t)$,
with certain positive constants $C', \gamma'$, 
which is independent of the initial value $a$.

RDEs (\ref{eq.rde}) and (\ref{eq.rde.der.sh}) combined as one has a linear growth coefficient.
So, it is not at all clear whether the solution $J$ exists or not.
(It is not very difficult to show the uniqueness, if a (global) solution exists.)
There are some papers (e.g., Lejay \cite{lej1}) which study Lyons' continuity theorem for RDEs with unbounded coefficients.
But, the case of the derivative equations (\ref{eq.rde})--(\ref{eq.rde.der.sh}) does not seem to be included.  

In solving an RDE with an unbounded coefficient,  
the most difficult part is always how to control the first level path of a solution.
In this case, however, 
thanks to the special form in (\ref{eq.rde.mat.sh}) and the series representation 
(\ref{eq.series}) below,  it is possible to prove existence of a unique solution $(Y, J)$.

Let us recall how this is solved in \cite{aida}:
\\
\\
\noindent
{\bf (Step 1)}
let us first consider the case when
$X$ is a smooth rough path lying above $x \in C^{1-var}_0 ([0,1], {\cal V})$.
Then, ODEs (\ref{eq.rde}) and (\ref{eq.rde.der.sh}) has a unique solution $t \mapsto (y_t, j_t)$
in 1-variational sense.
Moreover, it is well-known that $j_t$ can be written explicitly as follows;
\begin{equation}\label{eq.series}
{\rm Id} +j_t = ( {\rm Id} + j_s) 
\Bigl(
{\rm Id} +\sum_{k=1}^{\infty} A_{k; s,t}
\Bigr),
\qquad\qquad
(s \le t)
\end{equation}
where $A_k$ is given by
\begin{equation}\label{eq.def.mk}
A_{k; s,t}
=
\int_{s < t_1 < \cdots <t_k <t}  dM_{t_k} \cdots dM_{t_1},
\end{equation}
and $M_t$ on the right hand side
is given by (\ref{eq.rde.M})  in 1-variational sense
(with $X$ and $Y$ being  replaced with $x$ and $y$, resp.).
Notice the order of the product of $M_{t_j}$'s on the right hand side.
\\
\\
\noindent
{\bf (Step 2)}
Fortunately, $A_k$ is written in the form of interated integral.
So, the series representaion in (\ref{eq.series})--(\ref{eq.def.mk}) fits well with rough path theory.
The following argument is quite similar to the ``fundamental theorem 
of rough path theory''  (Theorem 3.1.2, \cite{lq}), 
which states that one can obtain the $i$th level path ($i \ge 3$) from the first 
and the second level paths.

Note also  that, in the same way as in the fundamental theorem,
the map 
$
x \mapsto M \mapsto j
$
extends continuously  with respect to the topology of $G\Omega_p ({\cal V})$.
\\
\\
\noindent
{\bf (Step 3)}
Suppose that $x \in C^{1-var}_0 ([0,1], {\cal V})$ satisfies that 
$\|X^1 \|_{p-var}^p + \|X^2\|_{p/2-var}^{p/2} \le R$ for $R>0$.
Then, it is shown that
$\sup_{0 \le t \le 1}|j_t| \le C_R <\infty$ for some positve constant $C_R$.

So, when we try to solve RDEs (\ref{eq.rde})--(\ref{eq.rde.der.sh}) 
for such $x$, the (local) solution $(Y,J)$ coincides with (the lift of) 
the solution $(y,j)$ in the usual sense 
and it (=its first level path) does not get out a large ball of radius $C'_{R} >0$.
It is also shown in \cite{aida} that,  for such an $x$ and its lift $X$,
Lipschits property for the map holds;
$$
X \mapsto (Y,J) \in G\Omega_p ( {\cal W} \oplus L({\cal W})).
$$
Therefore, this map naturally extends to a one from 
$\{X \in G\Omega_p ( {\cal V})  ~|~ \|X^1 \|_{p-var}^p + \|X^2\|_{p/2-var}^{p/2} \le R \}$ 
to $G\Omega_p ( {\cal W} \oplus L({\cal W}))$.
Since $R>0$ is arbitrary, this map is defined for any $X \in G\Omega_p ( {\cal V})$.
\\
\\
Summing up, we have the following proposition in \cite{aida}.
\begin{pr}\label{pr.sol.j} 
RDEs (\ref{eq.rde})--(\ref{eq.rde.der}) has a unique solution for any $X$
and the map 
\[
X\in  G\Omega_p ( {\cal V}) \mapsto (Y,J) \in G\Omega_p ( {\cal W} \oplus L({\cal W}))
\]
is (locally Lipschitz) continuous.
Moreover, $J^1$ i.e., the first level path of $J$, 
admits a series representation as in  (\ref{eq.series})--(\ref{eq.def.mk}).
\end{pr}


\section{Moment estimate of $J^1$}
In this section we prove that $L^r$-momoent of $J^1$ 
is finite when $X=W$,
by using the series in is finite by using the series representation 
in (\ref{eq.series})--(\ref{eq.def.mk}).

Let $M =(M^1,M^2)\in G\Omega_p (L({\cal W}))$.
Then, by the fundamental theorem of rough path theory,
we can construct $M^3, M^4, \ldots$.
When $M$ is a smooth rough path lying above $m$, 
then $M^k$ coincides with the iterated Stieltjes integral, i.e., 
\[
M^k_{s,t} =\int_{s<t_1 <\cdots < t_k<t }   dm_{t_1} \otimes \cdots \otimes dm_{t_k},
\qquad
\qquad
(k\ge 1).
\]
$A_k$ in (\ref{eq.def.mk}) is dominated by $M^k$. 
\begin{lm}\label{lm.iterA}
Set 
\[
A_{k: s,t} =\int_{s<t_1 <\cdots < t_k<t }   dm_{t_k}  \cdots  dm_{t_1}
\]
for a smooth rough path  $M =(M^1,M^2)\in G\Omega_p (L({\cal W}))$
lying above $m$.
Then, $M \mapsto A_{k: s,t}$ extends to a continuous map from $G\Omega_p (L({\cal W}))$
and the following inequality holds;
\[
|A_{k: s,t}| \le | M^k_{s,t} |  \qquad\qquad \mbox{for all $k$ and $s \le t$}
\]
\end{lm}

\Proof
Set $T: L({\cal W})^{\otimes k} \to L({\cal W})$ by 
$T(a_1 \otimes \cdots \otimes a_k )=a_k \cdots a_1$.
A basic propety of the projective tensor norm, the operator norm of $T$ is $1$.
Noting that $A_{k: s,t} =T( M^k_{s,t})$, we can easily prove the lemma.
\QED



It is well-known how fast $L^r$-norm of $W^i~(i=1,2)$ grows as $r \to \infty$,
from which one can obtain growth of $M^i~(i=1,2)$.
\begin{lm}\label{lm.fer}
There is a positive constant $C$ such that,
\begin{equation}\label{ineq.lrgrow}
{\mathbb E}[ \| W^i    \|_{i/p-Hld}^r   ]^{1/r}  \le C r^{i/2}     \qquad  \mbox{ for all $r \ge 1$ and $i=1,2$.} 
\end{equation}
\end{lm}

\Proof
Recall that a Fernique-type theorem holds for Brownian rough path;
there exists a positive constant $\beta$ such that ${\mathbb E} [ \exp( \beta (\| W^1   \|_{1/p-Hld} +\| W^i    \|_{i/p-Hld}^{1/2} )^2 )] <\infty$. 
(This type of theorem was first shown in \cite{lqz}. 
For similar results, see also \cite{ikldp, ms, fo} for instance.)
Dereich \cite{der} gave a proof under Assumption {\bf (A1)}.

We will show that, in general, if a random variable $Z \ge 0$ defined on a certain probability space
with  Fernique-type integrability condition
${\mathbb E} [ e^{\beta Z^2}] <\infty$,
then $\|  Z\|_{L^r} =O (\sqrt{r})$.
(In the following, we assume $\beta >1/2$ for simplicity. 
Otherwise, we take constant multiple of $Z$ instead of $Z$ itself.)
By Chebychev's inequality,  there exists a constant $C>0$ such that 
${\mathbb P} (Z \ge \eta) \le C e^{- \beta \eta^2}$ for all $\eta>0$.
Let $\ve >0$ be sufficiently small so that $\alpha := \beta -\ve >1/2$.

Then, for $r \ge 1$,
\begin{eqnarray*}
{\mathbb E} [Z^r]   &\le& \sum_{n=0}^{\infty}  (n+1)^r {\mathbb P} ( n < Z \le n+1 )
\nn\\
&\le&
C  \sum_{n=0}^{\infty}  (n+1)^r  e^{- \beta n^2}
\nn\\
&\le&
C r^{r/2}  \sum_{n=0}^{\infty}  \exp \bigl(  - \frac{r \log r}{2}    + r \log (n+1) -\alpha n^2  \bigr)   e^{-\ve n^2}.
\end{eqnarray*}
So, it suffices to show that $f(n,r):= - \frac{r \log r}{2}    + r \log (n+1) -\alpha n^2$ is bounded from above.
$\partial f / \partial n (n, r) = -2\alpha n+r /(n+1)$.
It is easy to see that, as a function of $n$, $f$ takes it maximum at $n = ( -1 + \sqrt{1 +2\alpha^{-1} r} )/2$.
Hence, 
\begin{equation}\label{ineq.fujiko}
f(n, r) 
\le 
 - \frac{r \log r}{2}   
 + r \log \Bigl( \frac{1 + \sqrt{1 + 2\alpha^{-1}  r}}{2} \Bigr)
- \alpha \Bigl(   \frac{ -1 + \sqrt{1+2\alpha^{-1}  r} }{2}   \Bigr)^2.
\end{equation}
It is easy to see that, as $r \to \infty$, 
\[
 \alpha \Bigl(   \frac{ -1 + \sqrt{1+2\alpha^{-1}  r} }{2}   \Bigr)^2 \approx \frac{r}{2}, 
 \quad
 \frac{1 + \sqrt{1 + 2\alpha^{-1}  r}}{2}  \approx \sqrt{ \frac{r}{2 \alpha}  }.
 \]
Take $\delta >0$ so small that $\log (1 +\delta) < 1/4$.
There exists $r_0 >0$ such that, 
for all  $r \ge r_0$, 
\[
\alpha \Bigl(   \frac{ -1 + \sqrt{1+2\alpha^{-1}  r} }{2}   \Bigr)^2 \ge  \frac{r}{4}, 
 \quad
 \log \Bigl( \frac{1 + \sqrt{1 + 2\alpha^{-1}  r}}{2} \Bigr) 
    \le    \log \Bigl(  (1+\delta) \sqrt{ \frac{r}{2 \alpha}  }  \Bigr).
    \]
  So, the right hand side of (\ref{ineq.fujiko}) is dominated by
  \[
  r \Bigl(   \log (1+\delta) - \frac12 \log (2\alpha)  -\frac 14 \Bigr)  \le 0,
  \qquad
  \mbox{ (for any $r \ge r_0$), }
  \]
  which implies that the right hand side of (\ref{ineq.fujiko}) is bounded from above.
\QED



\begin{lm}\label{lm.mom.M12}
When $X=W$ (i.e., Brownian rough path), 
there are positive constants $c$ and $\alpha$ independent of $r \ge 1$,
which satisfy that, for all $s \le t$, 
\begin{eqnarray*}
{\mathbb E} [|M^1_{s,t}|^r]^{1/r} &\le& \frac{ r^{\alpha} |c(t-s)|^{1/p}}{ \beta_p (1/p)!}
\\
{\mathbb E} [|M^2_{s,t}|^r]^{1/r} 
&\le& 
\frac{ r^{2\alpha} |c(t-s)|^{2/p}}{ \beta_p (2/p)!}.
\end{eqnarray*}
Here, $\beta_p$ is  a  positive constant such that
$$
\beta_p  \ge p^2 \Bigl(  1+ \sum_{j=3}^{\infty} \bigl( \frac{2}{j-2} \bigr)^{ 3/p} \Bigr).
$$ 
(Any such a positive constant which satisfies the above inequality  will do.)
\end{lm}

\Proof
We can choose 
$\omega(s,t)= ( \|W^1\|^p_{1/p-Hld} +  \|W^2\|^{p/2}_{2/p-Hld} )(t-s)$
as a control of $W$.
Then, as is explained in (\ref{eq.rde}), $M$  is controlled by $\omega '(s,t) =C'(1+\omega(0,1)^{\gamma '}) \omega(s,t)$,
where $C', \gamma'$ are positive constants.
Now, by using inequality (\ref{ineq.lrgrow}) in Lemma \ref{lm.fer} and choosing a suitable  constant  $c>0$, 
we can see the lemma holds for $2\alpha = 1 +( \gamma' /p)$.
\QED

\begin{pr}\label{lm.mom.Mall}
Let $X=W$ and set $\eta(s,t)= c(t-s)$, where $c$ is a positive constant 
as in Lemma \ref{lm.mom.M12} above.
Then, for all $k =1,2,3,\ldots$, $r \ge 1$, and $s \le t$, it holds that
\begin{eqnarray}\label{ineq.ind}
{\mathbb E} [|M^k_{s,t}|^r]^{1/r} \le
 \frac{ r^{k\alpha} \eta(s,t)^{k/p}}{ \beta_p (k/p)!}.
\end{eqnarray}
\end{pr}

\Proof
We use induction. The cases $k=1,2$ were already shown.
Assume the inequality (\ref{ineq.ind}) up to $k-1$ and let us prove (\ref{ineq.ind}) for $k$. 
%
%

Let us write
$M^j_{s,t}=M(a, X)^j_{s,t}$, where $X$ is the driving rough path 
and $a\in {\cal W}$ is the initial condition of the RDE. 
For a dyadic rational number
$u \in [0,1]$, we set $\tilde{w}_t = w_{u+t} -w_u$ for $0 \le t \le 1-u$.
Since $\tilde{w}$ and  $\{ w_s \}_{0 \le s \le u}$ are independent,
$\tilde{W}$, the lift of $\tilde{w}$, and $\{ W_{s,s'} \}_{0 \le s \le s' \le u}$
are independent.
For $u <t$,
$M(a, W)^j_{u,t}= M(a+ Y(a,W)^1_{0,u}, \tilde{W} )^j_{0,t-u}$.

From this independece, we see that, for $1 \le j \le k-1$,
\begin{eqnarray}
{\mathbb E}[ | M^j_{s,u} \otimes M^{k-j}_{u,t} |^r ]^{1/r}
&\le&
{\mathbb E}
\bigl[ | M(a,W)^j_{s,u}|^r \cdot  \tilde{{\mathbb E}} 
[|M (a+ Y(a,W)^1_{0,u}, \tilde{W} )^{k-j}_{0,t-u} |^r]  \bigr]^{1/r}
\nn
\\
&\le&
{\mathbb E}
\bigl[ | M(a,W)^j_{s,u}|^r   
\cdot
\Bigl( 
\frac{ r^{(k-j)\alpha} \eta(0,t-u)^{(k-j)/p}}{ \beta_p ((k-j)/p)!}
\Bigr)^r
\bigr]^{1/r}
\nn
\\
&\le&
 \frac{ r^{j\alpha} \eta(s,u)^{j/p}}{ \beta_p (j/p)!}
\cdot
\frac{ r^{(k-j)\alpha} \eta(u,t)^{(k-j)/p}}{ \beta_p ((k-j)/p)!}.
\label{ineq.lr.mm}
\end{eqnarray}
Here, $\tilde{{\mathbb E}}$ denotes the expectation with respect to $\tilde{W}$.
By taking limit, we can easily see that (\ref{ineq.lr.mm}) holds for any $s<u<t$.

Let ${\cal P}= \{ s =t_0 <t_1 < \cdots < t_L=t\}$ be a partition
of the interval $[s,t]$ and set 
\[
M^k_{s,t} ({\cal P}) 
:=\sum_{j=1}^{k-1} \sum_{i=1}^{L}
M^j_{s,t_{i-1}} \otimes M^{k-j}_{t_{i-1} ,t_i}
\] 
Recall that, by the fundamental theorem of rough path theory, $M^k_{s,t}$
is obtained as the limit;
$M^k_{s,t} =\lim_{|{\cal P}| \searrow 0} M^k_{s,t} ({\cal P}). $

It is well-known that there exists $t_l \in {\cal P}$ such that 
$\eta(t_{l-1},t_{l+1}) \le 2 \eta(s,t)/ (L-1)$ if $L \ge 3$.
(When $L=2$, we have a trivial (in)equality.)
By straight forward computation, 
\begin{equation}\label{eq.rem.one}
M^k_{s,t} ({\cal P}) - M^k_{s,t} ({\cal P} \setminus \{t_l \})
=
\sum_{j=1}^{k-1} M^{j}_{t_{l-1} ,t_l}\otimes M^{k-j}_{t_{l} ,t_{l+1}}.
\end{equation}

Now, we use the binomial inequality (also known as the neo-classical inequality),
which states that, for any $p \ge 1,  a, b \in [\infty), k=1,2,\ldots,$, 
\begin{equation}\label{ineq.binom}
\sum_{j=0}^k  \frac{a^{j/p}  b^{(k-j)/p}}
{\bigl(  \frac{j}{p}\bigr)! \bigl(  \frac{k-j}{p}\bigr)! }
\le
p^2 \frac{(a+b)^{k/p}}{\bigl(  \frac{k}{p}\bigr)! }.
\end{equation}
(By the way, Hara and Hino \cite{hh} recently proved that the best constant on the right hand side of
(\ref{ineq.binom}) is $p$, not $p^2$.)


By taking $L^r$-norm of (\ref{eq.rem.one}) when $X=W$, we see that
\begin{eqnarray}
\| M^k_{s,t} ({\cal P}) - M^k_{s,t} ({\cal P} \setminus \{t_l \}) \|_{L^r}
&\le&
\sum_{j=1}^{k-1} \| M^{j}_{t_{l-1} ,t_l}\otimes M^{k-j}_{t_{l} ,t_{l+1}} \|_{L^r}
\nn
\\
&\le&
\sum_{j=1}^{k-1}
\frac{ r^{j\alpha} \eta(t_{l-1} ,t_l  )^{j/p}}{ \beta_p (j/p)!}
\cdot
\frac{ r^{(k-j)\alpha} \eta( t_{l} ,t_{l+1}  )^{(k-j)/p} }  { \beta_p ((k-j)/p)!}
\nn
\\
&\le& 
\bigl( \frac{p}{\beta_p} \bigr)^2
\frac{r^{k\alpha} \eta(t_{l-1} ,t_{l+1})^{k/p} }{(k/p)!}
\nn
\\
&\le& 
\bigl( \frac{p}{\beta_p} \bigr)^2
\bigl( \frac{2}{L-1} \bigr)^{k/p}
\frac{r^{k\alpha} \eta(s ,t)^{k/p} }{(k/p)!}
\label{ineq.biai}
\end{eqnarray}
for $L \ge 3$.
Hence, we see from the definition of $\beta_p$ that
\begin{eqnarray}
\| M^k_{s,t} ({\cal P})    \|_{L^r}
&\le&
\bigl( \frac{p}{\beta_p} \bigr)^2
\Bigl( 1 +  \sum_{L=3}^{\infty}  \bigl( \frac{2}{L-1} \bigr)^{k/p}   \Bigr)
\frac{r^{k\alpha} \eta(s ,t)^{k/p} }{(k/p)!}
\nn
\\
&\le&
\bigl( \frac{p}{\beta_p} \bigr)^2
\Bigl( 1 +  \sum_{L=3}^{\infty}  \bigl( \frac{2}{L-1} \bigr)^{3/p}   \Bigr)
\frac{r^{k\alpha} \eta(s ,t)^{k/p} }{(k/p)!}
%
\le
\frac{r^{k\alpha} \eta(s ,t)^{k/p} }{ \beta_p  (k/p)!},
\nn
\end{eqnarray}
which implies
\begin{equation}\label{ineq.biai3}
\| M^k_{s,t}    \|_{L^r}  \le \frac{r^{k\alpha} \eta(s ,t)^{k/p} }{ \beta_p  (k/p)!} =  \frac{  \{r^{\alpha}  C^{1/p} \}^k  (t-s)^{k/p} }{ \beta_p  (k/p)!}.
\end{equation}
This completes the proof.
\QED

\begin{re}\label{re.fracok}
notice that, in the proof of above lemma, we used the independent increment property of Brownian motion.
So, if the driving process is replaced with fractional Brownian motion,  this proof fails.
\end{re}

\begin{pr}\label{pr.convA}
There exists a positive constant $C=C_r$ which is independent of $s, t$ such that
$
 \sum_{ k=1}^{\infty}  \| A_{k: s,t}    \|_{L^r}  \le  C (t-s)^{1/p}
$
for all $s<t$ and  
$
\| {\rm Id}_{{\cal W}}  + J^1_{0,t}    \|_{L^r}  \le  C
$
for all $t$.
\end{pr}

\Proof
Recall Stirling's formula;
\[
\lambda ! := \Gamma (\lambda +1)  \sim  \sqrt{2 \pi \lambda}  \lambda^{\lambda} e^{- \lambda}
 \qquad \mbox{  as $\lambda \nearrow \infty$.}
\]
Then, the first estimate holds for $C= \sum_k   \{r^{\alpha}  C^{1/p} \}^k /  \{\beta_p  (k/p)!  \}$. 
The second estimate is clear from the first one since 
${\rm Id}_{{\cal W}}  + J^1_{0,t}   =  {\rm Id}_{{\cal W}}  +  \sum_{ k=1}^{\infty}  A_{k: 0,t} $.
\QED

Set
\[
\|  \psi \|_{m , \theta} := \Bigl(  \iint_{0 < s < t < 1}  \frac{ \| \psi_t - \psi_s \|_{{\cal W}}^{m}}{|t-s|^{2+m \theta}}   dsdt  \Bigr)^{1/m},
\qquad
\quad
m \in [1, \infty), \theta \in (0,1],
\]
for a continious path $\psi$ in the usual sesne, 
which takes its values in a Banach space ${\cal W}$ and starts at $0$. 
It is well-known that this norm is stronger than the H\"older norm;
there exitsts a positive constant $C_{m, \theta}$ such that, for all $\psi$,
$
\|  \psi \|_{1/p-Hld} \le C_{m, \theta} \| \psi \|_{m , 1/p}.
$

Now we give the main result of this paper, 
which states that the first level path of the derivative equation (\ref{eq.rde.der}) is $L^r$-integrable for any $r \ge 1$,
provided $X=W$.
\begin{tm}\label{tm.main}
Assume {\bf (A1)} and  $\sigma : {\cal W} \to L({\cal V}, {\cal W})$ be $C^4_b$.
Consider RDEs  (\ref{eq.rde}) and  (\ref{eq.rde.der.sh})  
with the driving rough path 
being Brownian rough path,  i.e., $X=W$.
We denote by  $J^1$ the first level path of the solution of  the second RDE  (\ref{eq.rde.der.sh}).
Then, 
${\mathbb E} [ \|  J^1 \|^r_{1/p-Hld}  ]  <\infty$ for any $r \ge 1$.
\end{tm}

\Proof
From the series representation (\ref{eq.series}), it is clear that
$
J^1_{s,t} =  ( {\rm Id} + J^1_{0,t } )  \sum_{k=1}^{\infty} A_{k:s,t}.
$
 We obtain from Proposition \ref{pr.convA} that
 $$ \|  J^1_{s,t} \|_{L^r }  = \|  {\rm Id} + J^1_{0,t } \|_{L^{2r }}  \cdot \|   \sum_{k=1}^{\infty} A_{k:s,t}  \|_{L^{2r }}  \le C \sqrt{t-s}.$$ 
If $r >1$ is so large that  $2 + r  ( p^{-1} - 2^{-1} )  <1$, 
then we can easily see that 
\begin{eqnarray*} 
{\mathbb E} [  \|  J^1 \|^r_{1/p-Hld}  ]  &\le&  {\mathbb E} [  \|  J^1 \|^r_{r, 1/p}  ]  
\le
  \iint_{0 < s < t < 1}  \frac{  {\mathbb E}[  |J^1_{s,t}|^r ] } {|t-s|^{2+r /p}}   dsdt 
  \nn\\
  &\le&
  C   \iint_{0 < s < t < 1}  \frac{ 1} { |t-s|^{2+r (1/p  -1/2 )}}   dsdt  <\infty.
     \end{eqnarray*} 
 Thus, we have shown the theorem.
 \QED



\begin{thebibliography}{99}

\bibitem{aida}
Aida, S.; 
{\it T. Lyons no renzokusei teiri no shoumei ni tsuite,}
 (in Japanese:  "On the  proof of T. Lyons' contiuity theorem"), 
 unpublished.


\bibitem{ba}
Ben Arous, G.; 
Methods de Laplace et de la phase stationnaire sur l'espace de Wiener.  
 Stochastics  25  (1988),  no. 3, 125--153.

\bibitem{der}
Dereich, S.;
Rough paths analysis of general Banach space-valued Wiener processes. 
 J. Funct. Anal.  258  (2010),  no. 9, 2910--2936.


\bibitem{fo}
Friz, P.; Oberhauser, H.;
A generalized Fernique theorem and applications,
Proc. Amer. Math. Soc.  138  (2010), 3679-3688. 

\bibitem{hh}
Hara, K.; Hino. M.;
Fractional order Taylor's series and the neo-classical inequality,
Bull.  Lond. Math. Soc.  42  (2010),  467--477.


\bibitem{ikldp}
Inahama, Y.; Kawabi, H,;
Large deviations for heat kernel measures on loop spaces via rough paths,
J. London Math. Soc. (2)  73  (2006),  no. 3, 797--816. 

\bibitem{ik}
Inahama, Y.; Kawabi, H,;
Stationary  phase for solutions of rough differential equations.
In preparetion.


\bibitem{lqz}
Ledoux, M.; Qian, Z.; Zhang, T.;
Large deviations and support theorem for diffusion processes via rough paths,
Stochastic Process. Appl. 102 (2002), no. 2, 265--283. 

\bibitem{lej1}
Lejay, A.;
On rough differential equations.  Electron. J. Probab.  14  (2009), no. 12, 341--364.



\bibitem{lq}
Lyons, T.; Qian, Z.;
System control and rough paths. Oxford Mathematical Monographs. Oxford Science Publications. Oxford University Press, Oxford, 2002.


\bibitem{ms}
Millet, A.; Sanz-Sol\'e, M.;
Large deviations for rough paths of the fractional Brownian motion,
  Ann. Inst. H. Poincar\'e Probab. Statist.  42  (2006),  no. 2, 245--271.


\end{thebibliography}
\end{document}